\documentclass{amsart}
\usepackage{graphicx}

\newtheorem{theorem}{Theorem}
\newcommand{\bt}{\begin{theorem}}
\newcommand{\et}{\end{theorem}}

\newtheorem{lemma}{Lemma}
\newcommand{\bl}{\begin{lemma}}
\newcommand{\el}{\end{lemma}}

\newtheorem{corollary}{Corollary}
\newcommand{\bc}{\begin{corollary}}
\newcommand{\ec}{\end{corollary}}

\newtheorem*{conjecture}{Conjecture}
\newcommand{\bconj}{\begin{conjecture}}
\newcommand{\econj}{\end{conjecture}}

\newtheorem{example}{Example}
\newcommand{\bex}{\begin{example}}
\newcommand{\eex}{\end{example}}

\newtheorem{problem}{Problem} 
\newcommand{\bp}{\begin{problem}}
\newcommand{\ep}{\end{problem}}

\newcommand{\beq}{\begin{equation}}
\newcommand{\eeq}{\end{equation}}
\newcommand{\benum}{\begin{enumerate}}
\newcommand{\eenum}{\end{enumerate}}
\newcommand{\ba}{\begin{array}}
\newcommand{\ea}{\end{array}}
\newcommand{\Z}{\ensuremath{\mathbf Z}}

\newcommand{\PP}{\ensuremath{\mathbf P}}

\newcommand{\F}{\ensuremath{\mathbf F }}
\newcommand{\Fp}{\ensuremath{ {\mathbf F}_p }}

\begin{document}
\title{Heights on the finite projective line}
\author{Melvyn B. Nathanson}
\address{Department of Mathematics\\
Lehman College (CUNY)\\
Bronx, New York 10468}
\email{melvyn.nathanson@lehman.cuny.edu}
\thanks{The work of M.B.N. was supported in part by grants from the NSA Mathematical Sciences Program and the PSC-CUNY Research Award Program.}


\begin{abstract}
Define the height function $h(a) = \min\{k+(ka\mod p):k=1,2,\ldots,p-1\}$ for $a \in \{0,1,\ldots,p-1.\}$  It is proved that the height has peaks at $p, (p+1)/2,$ and $(p+c)/3,$ that these peaks occur at $a= [p/3], (p-3)/2, (p-1)/2, [2p/3], p-3,p-2,$ and $p-1,$ and that $h(a) \leq p/3$ for all other values of $a$.
\end{abstract}

\date{\today}
\maketitle

\section{Heights on finite projective spaces}
Let $p$ be an odd prime and let $\Fp=\Z/p\Z$ and  $\F_p^{\ast} = \F_p\setminus\{p\Z\}$.  For 
$d \geq 2$, we define an equivalence relation on the set of nonzero 
$d$-tuples in $\F_p^{d}$ 
as follows: 
$(a_1+p\Z,\ldots,a_d+p\Z) \sim (b_1+p\Z,\ldots,b_d+p\Z)$ 
if there exists $k\in \{1,2,\ldots,p-1\}$ such that 
$(b_1+p\Z,\ldots,b_{d}+p\Z)=(ka_1+p\Z,\ldots,ka_{d}+p\Z)$.
 We denote the equivalence class of $(a_1+p\Z,\ldots,a_{d}+p\Z)$
by $\langle a_1+p\Z,\ldots,a_{d}+p\Z\rangle$.
The set of equivalence classes is called the $(d-1)$-dimensional projective space over the field $\F_p$, and denoted $\PP^{d-1}(\F_p)$. 

For every integer $x$, we denote by  $x \mod p$  the least nonnegative  integer in the congruence class $x+p\Z$.  
We define the \emph{height} of the point  $\mathbf{a} = \langle a_1+p\Z,\ldots,a_{d}+p\Z\rangle \in \PP^{d-1}(\F_p)$ by
\[
h_p(\mathbf{a}) = \min \left\{\sum_{i=1}^d (ka_i \mod p) : k=1,\ldots,p-1\right\}.
\]
Nathanson and Sullivan~\cite{nath07p} introduced this definition of height in connection with a problem of Chudnovsky, Seymour, and Sullivan~\cite{chud-seym-sull07} in graph theory.

\begin{figure}  \label{HeightsLine:figure:4409}
\scalebox{.60}{\includegraphics{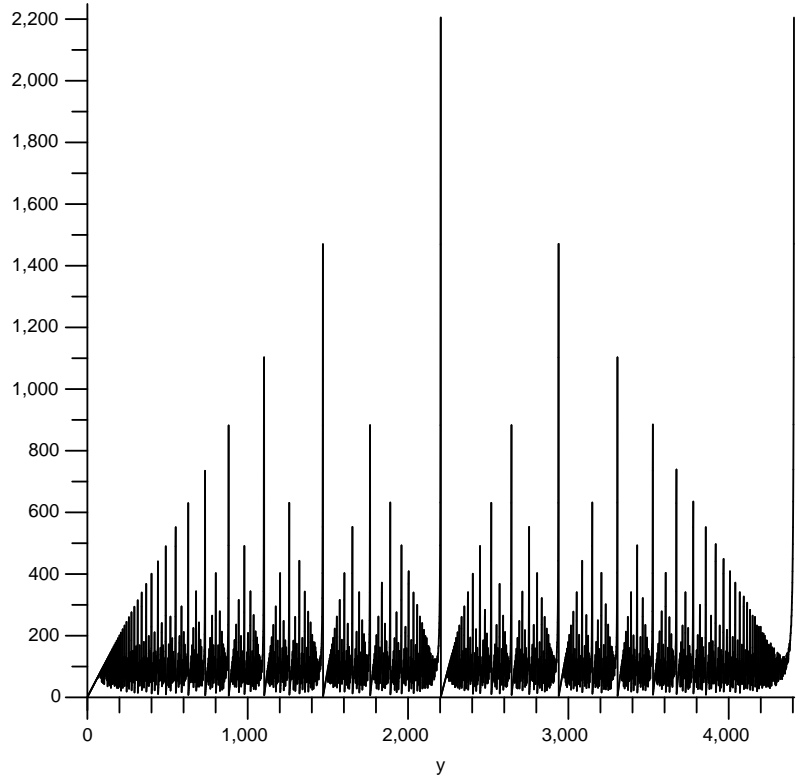}}
\caption{Heights for the prime 4409}
\end{figure}

The projective line $\PP^1(\F_p)$ consists of all equivalence classes of pairs $(a_1+p\Z,a_2+p\Z)$, where $a_1$ and $a_2$ are not both 0 modulo $p$.
If $a_1 \mod p =0$, then $\langle a_1+ p\Z,a_2+p\Z\rangle  = \langle p\Z,1+p\Z\rangle $ and $h_p(\langle p\Z,1+p\Z\rangle )=1.$
If $a_1\mod p\neq 0$, then there is an integer $a$ such that $a \equiv a_1^{-1}a_2 \pmod{p}$ and $\langle a_1+p\Z,a_2+p\Z\rangle  = \langle 1+p\Z,a+p\Z\rangle $.  
Thus, for all $\mathbf{a} \in \PP^1(\F_p)$, if $\mathbf{a} \neq \langle p\Z,1+ p\Z\rangle $, then $\mathbf{a} = \langle 1+p\Z,a+p\Z\rangle $ for some integer $a\in \{0,1,\ldots,p-1\}.$ 
For points on the projective line, we denote the height function by
\[
h(a) = h_p(\langle 1, a+p\Z \rangle).
\]
In Figure~\ref{HeightsLine:figure:4409} we graph the height function $h(a)$ for $a=0,1,\ldots, p-2$ for the  prime $p=4409.$   Note that $h(p-1)=p$ for all primes $p$.  The graph shows peaks at points $a \approx p/2, p/3, 2p/3,p/4,3p/4,\ldots$ with values $h(a) \approx p/2,p/3,p/3,p/4,p/4,\ldots,$ and suggests the conjecture that for every positive integer $b$ and all sufficiently large primes $p$, the range of the height function consists of integers that are equal approximately to $p,p/2,p/3,\ldots, p/b,$ together with a set of numbers less than $p/b$.  The goal of this paper is to prove this conjecture for $b=1,2,3.$

\section{Inequalities and identities}
\bl         \label{HeightsLine:lemma:line}
Let $p$ be an odd prime number and let $m$ and $\ell$  be integers such that $1 \leq \ell \leq m \leq p-1.$  If
\[
\frac{(\ell-1) p}{m} \leq a < \frac{\ell  p}{m}
\]
then
\[
h(a) \leq ma-( \ell - 1) p+m.
\]
 \el

\begin{proof}
Since
\[
0 \leq ma- (\ell -1) p < p,
\]
it follows that
\[
ma \mod p = ma - (\ell -1) p.
\]
The definition of height implies that 
\[
h( a  ) \leq m + (ma \mod p) = ma- (\ell -1)p+m.
\]
\end{proof}

\bl           \label{HeightsLine:lemma:hyperbola}
Let $p$ be an odd prime number and let $\ell$  and $m$ be integers such that  $1 \leq \ell \leq m \leq p-1.$  
Let
\[
H_{\ell,m,p}(x) =  \ell p - mx - 1 + \frac{mp}{\ell p - mx}
\]
for $x < \ell p/m.$  If
\[
\frac{(\ell-1) p}{m} < a < \frac{\ell  p -1}{m}
\]
then
\[
h(a) < H_{\ell,m,p}(a).
\]
Let $u$ and $v$ be real numbers such that 
\[
\frac{(\ell-1) p}{m} < u \leq \frac{\ell p}{m} - \sqrt{\frac{p}{m}} \leq v < \frac{\ell  p -1}{m}.
\]
If $u \leq a \leq v,$
then
\[
h(a) < \max( H_{\ell,m,p}(u), H_{\ell,m,p}(v)).
\]
\el

\begin{proof}
We have $1 <  \ell p - ma < p.$  
Dividing the prime number $p$ by $\ell p - ma$, we obtain
\[
p = q(\ell p - ma) + r
\]
where \[
1 \leq q = \left[ \frac{p}{\ell p - ma}\right] < \frac{p}{\ell p - ma} 
\]
and 
\[
1 \leq r \leq \ell p - ma -1 < p-1.
\]
Let $k = mq \pmod p.$ 
Then $1 \leq k \leq p-1$ and 
\[
m q a = p (\ell q -1) +r 
\]
and so
\[
ka \mod p = mq a  \mod p = r.
\]
It follows that
\begin{align*}
h(a) & \leq k + (ka \mod p) \\
& \leq mq + r  \\
& < \ell p - ma - 1 +  \frac{mp}{\ell p-ma} \\
& = H_{\ell,m,p}(a).
\end{align*}
The function $H_{\ell,m,p}(x)$ is decreasing for $x \leq \ell p /m - \sqrt{p/m}$ and increasing for $\ell p /m - \sqrt{p/m} \leq x < \ell p/m.$  It follows that if
\[
\frac{(\ell-1) p}{m} < u \leq  \frac{\ell p}{m} - \sqrt{\frac{p}{m}}  \leq v < \frac{\ell  p -1}{m}
\]
and $u \leq x \leq v,$ then 
$H_{\ell,m,p}(x) \leq \max(H_{\ell,m,p}(u), H_{\ell,m,p}(v)).$  
This completes the proof.
\end{proof}

\bl     \label{HeightsLine:lemma:hyperbola-ineq}
Let $p$ be an odd prime number and let $m$ and $\ell$  be integers such that  $1 \leq \ell \leq m \leq p-1.$  
If $m < w \leq \sqrt{mp}$ and
\[
\frac{\ell p}{m} - \frac{p}{w} \leq a \leq  \frac{\ell p}{m} - \frac{w}{m}
\]
then 
\[  
h(a) < \frac{mp}{w} + w-1.
\]
\el

\begin{proof}
We apply Lemma~\ref{HeightsLine:lemma:hyperbola}
with $u = \ell p/m - p/w$ and $v = (\ell p - w)/m.$ 
If $ m < w \leq \sqrt{mp},$ then 
\[
\frac{(\ell-1) p}{m} < 
\frac{\ell p}{m} - \frac{p}{w} 
 \leq  \frac{\ell p}{m} - \sqrt{\frac{p}{m}} \leq  
 \frac{\ell p}{m} - \frac{w}{m} < \frac{\ell  p -1}{m}.
\]
Since the function $H_{\ell,m,p}(w)$ satisfies the functional equation
\[
H_{\ell,m,p}\left( \frac{\ell p}{m} - \frac{p}{w} \right) 
= H_{\ell,m,p}\left( \frac{\ell p}{m} - \frac{w}{m} \right) 
= \frac{mp}{w}+w-1
\]
it follows that if
\[
\frac{\ell p}{m}-\frac{p}{w} \leq a \leq \frac{\ell p}{m} - \frac{w}{m}
\]
then
\[
h(a) < \max\left( H_{\ell,m,p}(u),  H_{\ell,m,p}(v) \right)  = \frac{mp}{w}+w-1.
\]
This completes the proof.
\end{proof}

\bl            \label{HeightsLine:lemma:p-b}
Let $p$ be an odd prime number and $b$ a positive integer.  If $p>(b-1)^2,$  then 
\[
h(p-b) = \frac{p}{b} + \frac{r(b-1)}{b}
\]
where $p\mod b = r.$
\el

\begin{proof}
For each positive integer $i$, divide $ip$ by $b$ to obtain
\[
ip = q_ib+r_i  \qquad\text{ with $0 \leq r_i \leq b-1.$}
\]
For each $i$, consider the set of integers $k$ that satisfy the inequality 
\[
\frac{(i-1)p}{b} < k \leq \frac{ip}{b} = q_i + \frac{r_i}{b} .
\]
Then $k \leq q_i.$  Every integer $k = 1,2,\ldots, p-1$ satisfies exactly one of these inequalities for some $i \leq b.$  If $k$ satisfies the $i$th inequality, then
\[
0 \leq ip-kb < p
\]
and 
\[
k(p-b) \mod p = ip-kb.
\]
It follows that 
\begin{align*}
k + (k(p-b) \mod p)  & = ip - k(b-1) \\
& \geq ip - q_i(b-1) \\
& = r_i + q_i \\
& = r_i + \frac{ip - r_i}{b} \\
& = \frac{ip}{b} + \frac{r_i(b-1)}{b} .
\end{align*}
The inequality $p > (b-1)^2$ implies that $(i+1)p+r_{i+1}(b-1) > ip + r_i (b-1)$ for $i \geq 1,$ and so
\[
h(p-b) = \min\left\{k(p-b) \mod p + k  : k=1,\ldots, p-1\right\} \geq \frac{p}{b} + \frac{r_1(b-1)}{b}.
\]
Choosing $i=1$ and $k=q_1$ gives $h(p-b ) = (p + r_1(b-1))/b.$
\end{proof}

\bl    \label{HeightsLine:lemma:(p-b)/2}
Let $p$ be an odd prime number.  Then
\[
h\left( \frac{p-1}{2}  \right) = \frac{p+1}{2}.
\]
If $b$ is an odd integer, $b \geq 3,$ and $p>(b-1)^2,$  then 
\[
h\left(\frac{p-b}{2}  \right) = \frac{p}{b} + \frac{(b-2)r}{2b}
\]
where $p+b \mod 2b = r.$
\el

\begin{proof}
Let $a = (p-1)/2$ and $k \in \{1,\ldots, p-1\}.$  If $k=2\ell$ for $\ell \in \{1,\ldots, (p-1)/2\},$ then
\[
ka \mod p = \ell(p-1)\mod p = p-\ell
\]
and 
\[
k + (ka\mod p) = p+\ell \geq p+1.
\]
If $k=2\ell-1$ for $\ell \in \{1,\ldots, (p-1)/2\},$ then
\[
ka \mod p = \left(\ell(p-1)- \frac{p-1}{2}\right)\mod p = \frac{p+1}{2}-\ell
\]
and 
\[
k + (ka\mod p) = \frac{p-1}{2}+\ell \geq \frac{p+1}{2}.
\]
It follows that $h((p-1)/2)=(p+1)/2.$

Let $b \geq 3$ and $k \in \{1,2,\ldots,p-1\}.$   If $k$ is even, then $k = 2\ell$ for $\ell \in \{ 1,\ldots, (p-1)/2\}.$  
For each positive integer $i$, divide $ip$ by $b$ to obtain
\[
ip = q_{0,i}b+r_{0,i}  \qquad\text{ with $0 \leq r_{0,i} \leq b-1.$}
\]
For each $i$, consider the set of integers $\ell$ that satisfy the inequality 
\[
\frac{(i-1)p}{b} < \ell \leq \frac{ip}{b} = q_{0,i} + \frac{r_{0,i}}{b} .
\]
Then $\ell \leq q_{0,i}.$  Every integer $\ell = 1,2,\ldots, (p-1)/2$ satisfies exactly one of these inequalities for some $i \leq (b+1)/2.$  If $\ell$ satisfies the $i$th inequality, then
\[
0 \leq ip-\ell b < p
\]
and 
\[
k\left( \frac{p-b}{2} \right) \mod p = \ell (p-b) \mod p = i p- \ell b.
\]
It follows that 
\begin{align*}
k + \left( k\left( \frac{p-b}{2} \right) \mod p \right) & = ip - \ell(b-2) \\
& \geq ip - q_{0,i}(b-2) \\
& = r_{0,i} + 2q_{0,i} \\
& = r_{0,i} + \frac{2ip - 2r_{0,i}}{b} \\
& = \frac{2ip}{b} + \frac{(b-2) r_{0,i}}{b} .
\end{align*}
The inequality $p > (b-1)^2$ implies that $2(i+1)p+ (b-2) r_{0,i+1} > 2ip + (b-2) r_{0,i} $ for $i \geq 1,$ and so, for even integers $k$, we have 
\beq    \label{HeightsLine:evencase}
\min\left\{k + \left(k(p-b) \mod p \right)  : k=2,4,\ldots, p-1\right\} \geq \frac{2p}{b} + \frac{(b-2) r_{0,1}}{b}.
\eeq

If $k$ is odd, then $k = 2\ell -1$ for some $\ell \in \{1,\ldots, (p-1)/2\}.$  
For each positive integer $i$, divide $(2i-1)p+b$ by $2b$ to obtain
\[
(2i-1)p+b = 2bq_{1,i}+r_{1,i}  \qquad\text{ with $0 \leq r_{1,i} \leq 2b-1.$}
\]
For each $i$, consider the set of integers $\ell$ that satisfy the inequality 
\[
0 \leq ip - \left( \ell b + \frac{p-b}{2}  \right)
< p,
\]
Then $\ell \leq q_{1,i}.$  Every integer $\ell = 1,2,\ldots, (p-1)/2$ satisfies exactly one of these inequalities for some $i \leq (b+1)/2.$  If $\ell$ satisfies the $i$th inequality, then
\[
(i-1)p < \ell b + \frac{p-b}{2} \leq ip
\]
and 
\[
k\left( \frac{p-b}{2} \right) \mod p 
=  \ell p - \left(  \ell b + \frac{p-b}{2} \right) \mod p 
= i p- \ell b - \frac{p-b}{2}.
\]
It follows that 
\begin{align*}
k + \left( k\left( \frac{p-b}{2} \right) \mod p \right) & = ip - \ell(b-2)- \frac{p-b}{2} - 1\\
& \geq ip - q_{1,i}(b-2) - \frac{p-b}{2} - 1 \\
& =  \frac{(2i-1)p}{2} - \left(\frac{  (2i-1)p +b-r_{1,i}   }{2b}  \right) (b-2) + \frac{b}{2} - 1 \\
& = \frac{(2i-1)p}{b} + \frac{(b-2)r_{1,i}}{2b} .
\end{align*}
The inequality $p > (b-1)^2$ implies that 
\[
2(2i+1)p+ (b-2)r_{1,i+1} > 2(2i-1)p + (b-2) r_{1,i}
\]
for $i \geq 1,$ and so, for odd integers $k$, we have 
\[
\min\left\{k + \left( k\left( \frac{p-b}{2} \right) \mod p \right)  : k=1,3,\ldots, p-2\right\} \geq \frac{p}{b} + \frac{(b-2)r_{1,1}}{2b}.
\]
Choosing $i=1$ and $\ell=q_{1,1},$ we obtain
\beq    \label{HeightsLine:oddcase}
\min\left\{k + \left( k\left( \frac{p-b}{2} \right) \mod p \right)  : k=1,3,\ldots, p-2\right\}
 = \frac{p}{b} + \frac{(b-2)r_{1,1}}{2b}.
\eeq
Observing that 
\[
\frac{2p}{b} + \frac{r_{0,1}(b-2)}{b} > \frac{p}{b} + \frac{(b-2)r_{1,1}}{2b}
\]
and comparing~\eqref{HeightsLine:evencase} and~\eqref{HeightsLine:oddcase}, we conclude that $h\left((p-b)/2 \right) = (p + (b-1)r_{1,1})/b.$
\end{proof}

\section{Peaks with height $p/b$ for $b = 1,2,3$}

\bt
Let $p \geq 17$ be an odd prime number, and let
\[
A_1 = \{p-1\}
\]
and
\[
A_2=\left\{ \left[\frac{p}{2}\right] , p-2\right\}.
\]  
Then $h(a)=p$ if and only if $a\in A_1$ and $h(a)=(p+1)/2$ if and only if $a \in A_2.$  
If
\[
a \in \{0,1,2,\ldots,p-1\} \setminus (A_1 \cup A_2)
\]
then
\[
h\left(a \right)  \leq \frac{p}{2}  
\]
\et

\begin{proof}
By Lemma~\ref{HeightsLine:lemma:line} with $m=\ell = 1,$ we have
\[
h\left( a \right) \leq a+1 \leq \frac{p-1}{2}
\]
for $0 \leq a \leq (p-3)/2.$  

By Lemma~\ref{HeightsLine:lemma:line} with $m=\ell = 2,$ we have
\[
h\left( a  \right) \leq 2a-p+2  \leq \frac{p-1}{2}
\]
for $(p+1)/2 \leq a  \leq (3p-5)/4.$  

By Lemma~\ref{HeightsLine:lemma:hyperbola-ineq} with $m=\ell = 1$ and $w=3,$ we have
\[
h(a) < \frac{p}{3}+2 \leq \frac{p+1}{2}
\]
for $2p/3 \leq a \leq p-3.$  Note that $2p/3 \leq (3p-5)/4$ for $p \geq 17.$  This completes the proof.
\end{proof}

\bt  \label{HeightsLine:theorem:p/3}
Let $p$ be a prime number, and let 
\[
A_3 = \left\{ \left[\frac{p}{3}\right], \frac{p-3}{2},  \left[\frac{2p}{3}\right],  p-3 \right\}.
\]
If $p \equiv 1 \pmod{3}$, then $h(a ) = (p+2)/3$ for $a \in A_3.$ 
If $p \equiv 2 \pmod{3}$, then $h(a) = (p+1)/3$ for $a = (p-2)/3$ or $(p-3)/2$ and $h(a) = (p+4)/3$ for $a = (2p-1)/3$ or $p-3.$  
If 
\[
a \in \{0,1,2,\ldots, p-1\} \setminus \left(A_1 \cup A_2 \cup A_3 \right)
\]
then
\[
h(a) \leq \frac{p}{3}.
\]
\et

\begin{proof}
Consider primes  $p \equiv 1 \pmod{3}$.  
Then $p+3 \equiv 4\pmod{6}$  and Lemma~\ref{HeightsLine:lemma:(p-b)/2} implies that  $h((p-3)/2) = (p+2)/3.$
By Lemma~\ref{HeightsLine:lemma:p-b}, $h(p-3) = (p+2)/3.$

We define $q = (p-1)/3.$  If $k\in \{1,\ldots,p-1\},$ then $k= 3\ell - r$ for $\ell \in \{1,\ldots,q\}$ and $r\in \{0,1,2\}.$  
We shall prove that if $a = q$ or $2q,$ then $h(a) = \min\{ k + (ka\mod p) : k=1,\ldots, p-1\} = q+1.$  Since $k + (ka\mod p) \geq k+1,$ it suffices to consider only $k \leq q$ or, equivalently, $\ell \leq (q+r)/3.$   We have 
\[
kq = \frac{(3\ell-r)(p-1)}{3} = \ell(p-1) -rq =\ell p -\ell -rq.
\]
Since
\[
0 \leq \ell + rq \leq \frac{(3r+1)q+r}{3} \leq \frac{7q+2}{3} < p
\]
it follows that $kq \mod p = p -\ell -rq$ and 
\[
k + (kq\mod p) = p+2\ell -r(q+1).
\]
This minimum value of this expression is $q+1,$ and occurs at $\ell = 1, r=2,$ and $k=1.$ 

Similarly,
\[
k(2q) = \frac{(3\ell-r)(2p-2)}{3} = \ell(2p-2) - 2rq = 2\ell p - 2\ell - 2rq.
\]
If $r=2,$ then $p < 2\ell+4q < 2p$ and 
$k(2q) \mod p = 2p-2\ell-4q.$  Then $k + (k(2q)\mod p) = 2p+\ell - 4q +2
= 2q+\ell+4$.
If $r=0$ or 1, then $2\ell + 2rq \leq 2\ell+2q \leq (8q+2)/3 < p$ and
$k(2q)\mod p = p-2\ell-2rq.$  Then $k + (k(2q)\mod p) = p+\ell-2rq-r$ 
has a minimum value of $q+1$ at $\ell=1,r=1$, and $k=2.$  Thus, $h(a)=q+1$ if $a \in A_3.$

Next we consider primes $p \equiv 2 \pmod{3}$.  
Then $p+3 \equiv 2 \pmod{6}$  and Lemma~\ref{HeightsLine:lemma:(p-b)/2} implies that  $h((p-3)/2) = (p+1)/3.$
By Lemma~\ref{HeightsLine:lemma:p-b}, $h(p-3) = (p+4)/3.$

Let $q=(p-2)/3.$  We shall prove that $h(q) = q+1.$  It suffices to consider only $k + (kq\mod p)$ for $k = 1,\ldots, q.$  Let $k = 3\ell - r$ for $r\in \{0,1,2\}.$  If $k\leq q,$ then $\ell \leq (q+r)/3.$  We have
\[
kq = \frac{(3\ell-r)(p-2)}{3} = \ell p - 2\ell -rq
\]
and
\[
0 < 2\ell + rq \leq \frac{(2+3r)q+2r}{3} \leq \frac{8q+4}{3} < p.
\]
It follows that $kq\mod p = p-2\ell-rq$ and so
\[
k+(kq\mod p) = p+\ell-rq-2 = (3-r)q+\ell.
\]
The minimum value of this expression is $q+1,$ and occurs when $\ell=1, r=2,$ and $k=1.$

We shall prove that $h(2q+1)  = q+2.$  Let $k = 3\ell - r \leq q+1,$ where $r \in \{0,1,2\}$ and $1 \leq \ell \leq (q+r+1)/3.$  Then
\[
k(2q+1) = \frac{(3\ell-r)(2p-1)}{3} = 2\ell p - \ell -2rq-r.
\]
If $r=2,$ then $p < \ell + 4q + 2 = p+q+\ell < 2p$ implies that
\[
k(2q+1) \mod  p = 2p -  ( p+q+\ell ) = p - q - \ell
\]
and 
\[
k + k(2q+1) \mod  p = p - q +2 \ell - 2 = 2q+2\ell \geq 2q+2.
\]
If $r=0$ or 1, then 
\[
0 < \ell + 2rq + r \leq \ell + 2q + 2 < p
\]
and
\[
k + (k(2q+1)\mod p) = (3\ell - r) + (p-\ell -2rq-r) = p+2\ell-2r(q+1).
\]
The minimum value of this expression is $q+2$ and occurs when $\ell = r = 1$ and $k=2.$

We have proved that if $a \in A_3,$ then $h(a) \in \{(p+1)/3,(p+2)/3,(p+4)/3\}.$  We shall prove that if $a \notin A_1 \cup A_2 \cup A_3,$ then $h(a) \leq p/3.$


Applying Lemma~\ref{HeightsLine:lemma:line} with $\ell = m=1,$
we obtain $h(a) \leq p/3$ for
\[
0\leq a \leq \frac{p-3}{3}.
\]
Applying Lemma~\ref{HeightsLine:lemma:line} with $\ell = 2$ and $m=3,$
we obtain $h(a) \leq p/3$ for
\[
\frac{p}{3}\leq a \leq \frac{4p-9}{9}.
\]
Applying Lemma~\ref{HeightsLine:lemma:hyperbola-ineq} with $\ell = 1, m = 2,$ and $w=7,$ we have 
\[
h(a) < \frac{2p}{7}+6 
\]
for 
\[
\frac{5p}{14} \leq a \leq \frac{p-7}{2}.
\]
By Lemma~\ref{HeightsLine:lemma:(p-b)/2}, we have
\[
h\left( \frac{p-5}{2}\right) \leq \frac{p+12}{5}.
\]
Since
\[
\frac{5p}{14} \leq \frac{4p-9}{9},
\]
it follows that if $a\leq a \leq (p-5)/2$ and $h(a)>p/3,$ then $a = [p/3].$

Applying Lemma~\ref{HeightsLine:lemma:line} with $\ell = m=2,$
we obtain $h(a) \leq p/3$ for
\[
\frac{p}{2} \leq a \leq \frac{2p-4}{3}.
\]
Applying Lemma~\ref{HeightsLine:lemma:line} with $\ell = m=3,$
we obtain $h(a) \leq p/3$ for
\[
\frac{2p}{3} \leq a \leq \frac{7p-9}{9}.
\]
Applying Lemma~\ref{HeightsLine:lemma:hyperbola-ineq} with $\ell = m = 1$ and $w=4,$ we have 
\[
h(a) < \frac{p}{4}+3 
\]
for 
\[
\frac{3p}{4} \leq a \leq p-4.
\]
Since
\[
\frac{3p}{4} \leq \frac{7p-9}{9}
\]
it follows that if $p/2\leq a \leq p-4$ and $h(a)>p/3,$ then $a = [2p/3].$
This completes the proof.
\end{proof}

\section{Problems}
It is an open problem  to describe the spectrum of the height function on projective spaces $\PP^{d-1}(\F_p)$ for $d \geq 3.$  It would also be interesting to have a reasonable definition of the height function for points in projective space over arbitrary finite fields.

\def\cprime{$'$} \def\cprime{$'$} \def\cprime{$'$}
\providecommand{\bysame}{\leavevmode\hbox to3em{\hrulefill}\thinspace}
\providecommand{\MR}{\relax\ifhmode\unskip\space\fi MR }
\providecommand{\MRhref}[2]{%
  \href{http://www.ams.org/mathscinet-getitem?mr=#1}{#2}
}
\providecommand{\href}[2]{#2}


\begin{thebibliography}{1}

\bibitem{chud-seym-sull07}
M.~Chudnovsky, P.~Seymour, and B.~D.~Sullivan, 
\emph{Cycles in dense digraphs},
  arXiv:math.CO/0702147, 2007.

\bibitem{nath07p}
M.~B. Nathanson and B.~D.~Sullivan,
\emph{Heights in finite projective space, and a problem on finite directed graphs}, 
arXiv: math.NT/0703418, 2007.

\end{thebibliography}
\end{document}